\documentclass[a4paper,notitlepage]{article}%
\usepackage{amssymb}
\usepackage{graphicx}
\usepackage{amsmath}
\usepackage{amsthm}
\usepackage{amsfonts}%
\providecommand{\U}[1]{\protect\rule{.1in}{.1in}}
\theoremstyle{plain}
\newtheorem{theorem}{Theorem}[section]

\newtheorem{lemma}[theorem]{Lemma}

\theoremstyle{definition}

\numberwithin{equation}{section}
\numberwithin{theorem}{section}
\ifx\pdfoutput\relax\let\pdfoutput=\undefined\fi
\newcount\msipdfoutput
\ifx\pdfoutput\undefined\else
\ifcase\pdfoutput\else
\ifx\paperwidth\undefined\else
\ifdim\paperheight=0pt\relax\else\pdfpageheight\paperheight\fi
\ifdim\paperwidth=0pt\relax\else\pdfpagewidth\paperwidth\fi
\fi\fi\fi
\begin{document}

\title{A local and nonlocal coupling model involving the $p$-Laplacian\thanks{2020
\textit{Mathematics Subject Classification}. 35R11, 45K05, 47G20.}
\thanks{\textit{Key words and phrases}. Local equations, nonlocal equations,
$p$-Laplacian} \thanks{Partially supported by Secyt-UNC grant
CB33620180100016} }
\author{Uriel Kaufmann\thanks{FaMAF, CIEM, Universidad Nacional de C\'{o}rdoba, (5000)
C\'{o}rdoba, Argentina. \textit{E-mail address: }uriel.kaufmann@unc.edu.ar},
Ra\'{u}l Vidal\thanks{Corresponding author. FaMAF, CIEM, Universidad Nacional
de C\'{o}rdoba, (5000) C\'{o}rdoba, Argentina. \textit{E-mail address:
}raul.vidal@unc.edu.ar}
\and \noindent}
\maketitle

\begin{abstract}
In this paper we extend some results presented in \cite{julio} to the case of
the $p$-Laplacian operator. More precisely, we consider a model that couples a
local $p$-Laplacian operator with a nonlocal $p$-Laplacian operator through
source terms in the equation. The resulting problem is associated with an
energy functional. We establish the existence and uniqueness of a solution,
which is obtained via the direct minimization of the corresponding energy functional.

\end{abstract}

\section{Introduction and main results}

Nonlocal models can describe phenomena that are not well represented by
classical PDE's, for example, problems which have long-range interactions
and/or discontinuities. For instance, in the context of diffusion, long-range
interactions effectively describe anomalous diffusion, while in the context of
mechanics, cracks formation results in material discontinuities.

Nonlocal operators are defined through integration against an appropriate
kernel, which implies that their values at a given point depend on the entire
domain rather than just a neighbourhood around that point, as is typical for
differential operators. One of the most important examples is the fractional Laplacian.

For general references on nonlocal models we refer e.g. to
\cite{BCRR15,CR23,CCR06,CERW07,CERW08,CD07,MD14,SL10,W02,Z04} and its
references, while the articles \cite{CY,ERS,GT,IR,IRA,MRT,SLL} focus on the
study of nonlocal $p$-Laplacian operators.

In recent years there has been growing interest in models that combine local
and nonlocal effects, as they are capable of capturing more complex and
realistic phenomena. In such cases, nonlocal effects may arise in certain
regions of the domain, while in other regions the behavior is governed by
classical differential operators. See, for instance,
\cite{ABR21,ABR23,DB,DLLT,GW,GQR,K15,SOR,SS} and the references therein.

The study of nonlinear partial differential equations with $p$-Laplacian
operators has gained significant attention due to its broad range of
applications in fields such as physics, engineering, and image processing. In
this work, we analyze an elliptic equation that couples the local
$p$-Laplacian operator with a nonlocal $p$-Laplacian operator through source
terms. This coupling results in a variational structure, and we establish the
existence and uniqueness of solutions by minimizing the corresponding energy
functional. Our results extend some of those presented in \cite{julio}, where
the classical Laplacian is considered both in its local and nonlocal forms.

For coupling local and nonlocal models the previous strategies treat the
coupling condition as an optimization objective (the goal is to minimize the
mismatch of the local and nonlocal solutions on the overlap of their
sub-domains). Another approach is based on a partitioned procedure as a
general coupling strategy for heterogeneous systems, the system is divided
into sub-problems in their respective sub-domains, which communicate with each
other via the transmission conditions. As far as we are aware, the literature
lacks studies addressing this approach for models that involve $p$-Laplacian operators.

\subsection{Statement of the main result}

We assume that $\Omega\subset\mathbb{R}^{N}$ is an open bounded domain, such
that $\Omega$ is divided into two disjoint subdomains: a local region that we
will denote by $\Omega_{\ell}$ and a nonlocal region, $\Omega_{n\ell}$. Thus
we have $\Omega_{\ell},\Omega_{n\ell}\subset\Omega\subset\mathbb{R}^{N}$ with
$\Omega=(\overline{\Omega}_{\ell}\cup\overline{\Omega}_{n\ell})^{\circ}$.
Further, we assume that:

(1) $\Omega_{\ell}$ is connected and has a Lipschitz boundary.

(2) $\Omega_{n\ell}$ is $\delta$-connected. As in \cite{julio}, for $\delta
>0$, we say that an open set $U\subset\mathbb{R}^{N}$ is $\delta$-connected if
it cannot be written as a disjoint union of two (relatively) open nontrivial
sets that are at distance greater or equal than $\delta$.

(3) $dist(\Omega_{\ell},\Omega_{n\ell})<\delta.$

Our aim is to consider the following local-nonlocal problem, under suitable
hypothesis on the nonlinearity $f$ and the kernel $J$:%
\begin{equation}
\left\{
\begin{array}
[c]{ll}%
f\left(  x,u\right)  =-\Delta_{p}u+%
{\textstyle\int\limits_{{\small \Omega}_{n\ell}}}
J(x,y)\left\vert u(x)-u(y)\right\vert ^{p-2}\left(  u(x)-u(y)\right)  dy &
\text{in }\Omega_{\ell},\\
\partial_{\nu}u=0 & \text{in }\partial\Omega_{\ell}\cap\Omega,\\
u=0 & \text{in }\partial\Omega\cap\partial\Omega_{\ell},
\end{array}
\right.  \label{ec1}%
\end{equation}
and the following nonlocal equation in $\Omega_{n\ell}$,%
\begin{equation}
\left\{
\begin{array}
[c]{ll}%
f\left(  x,u\right)  =\int_{\mathbb{R}^{N}\setminus\Omega_{n\ell}%
}J(x,y)\left\vert u(x)-u(y)\right\vert ^{p-2}\left(  u(x)-u(y)\right)
dy\smallskip & \text{in }\Omega_{n\ell},\\
\qquad\qquad\,\,+2\int_{\Omega_{n\ell}}J(x,y)\left\vert u(x)-u(y)\right\vert
^{p-2}\left(  u(x)-u(y)\right)  dy & \\
u=0 & \text{in }\mathbb{R}^{N}\setminus\Omega.
\end{array}
\right.  \label{ec2}%
\end{equation}

\noindent Here $f:\Omega\times\mathbb{R}\rightarrow\mathbb{R}$ is a
Carath\'{e}odory function (that is, $f\left(  x,\cdot\right)  $ is continuous
for $a.e.$ $x\in\Omega$ and $f\left(  \cdot,\xi\right)  $ is measurable for
all $\xi\in\mathbb{R}$) that satisfies the following growth condition:%
\[
\left\vert f\left(  x,\xi\right)  \right\vert \leq a\left(  x\right)
+b\left(  x\right)  \left\vert \xi\right\vert ^{q}\quad\text{for }a.e.\text{
}x\in\Omega\text{, }\xi\in\mathbb{R},\leqno{(f)}
\]
where $0\leq q<p-1$, and $a$ and $b$ are nonnegative functions such that $a\in
L^{s}(\Omega)$, with $s>p^{\prime}$ (where, as usual, $1/p+1/p^{\prime}=1$)
and $b\in L^{\infty}(\Omega)$. Regarding the hypothesis on $J$, we shall
assume that:

(J1) $J$ is symmetric, and there exists $C>0$ such that $J(z)>C$ for all $z$
such that $\left\vert z\right\vert \leq2\delta$.

(J2) Let $1<q<\infty$. For $u\in L^{q}$ we have that%
\[
T_{J,q}(u):=\int_{\Omega}J(x,y)u(y)dy,
\]
defines a compact operator in $L^{q}(\Omega_{n\ell})$. For sufficient
conditions on $J$ for $T_{J,q}$ to be a compact operator we refer e.g. to
\cite[Theorem 1]{GS} or \cite[Chapter VI]{brezis}.

Let $p>1$. We next consider the space
\[
W:=\left\{  u\in L^{p}\left(  \Omega\right)  :u_{\left\vert \Omega_{\ell
}\right.  }\in W^{1,p}\left(  \Omega_{\ell}\right)  ,\text{ }u=0\text{ in
}\mathbb{R}^{N}\setminus\Omega\right\}  ,
\]
which is a Banach space equipped with the norm
\[
\Vert u\Vert_{W}:=\Vert u\Vert_{L^{p}(\Omega)}+\Vert|\nabla u|\Vert
_{L^{p}(\Omega_{\ell})}.
\]

\noindent Let $\widetilde{F}\left(  x,u\right)  :=\int_{0}^{u}f\left(
x,\xi\right)  d\xi$. Defined in $W$ we have the energy functional
$E:W\rightarrow\mathbb{R}$ given by
\[
E(u):=\int_{\Omega_{\ell}}\frac{\left\vert \nabla u\right\vert ^{p}}{p}%
+\frac{1}{p}\int_{\Omega_{n\ell}}\int_{\mathbb{R}^{N}}J(x,y)\left\vert
u(x)-u(y)\right\vert ^{p}dydx-\int_{\Omega}\widetilde{F}\left(  x,u\right)
\,dx.
\]
It is easy to check that this functional is Fr\'{e}chet
differentiable.\smallskip

We can now state our main result:

\begin{theorem}
\label{T1} Let $p>1$, and assume (1), (2), (3), (f), (J1) and (J2). Then there
exists a minimizer of $E$ in $W$. Moreover, the minimizer is a weak solution
of (\ref{ec1}) and (\ref{ec2}). Furthermore, if $\xi\rightarrow f\left(
x,\xi\right)  $ is strictly concave for $a.e.$ $x\in\Omega$, then the
minimizer of the functional $E$ is unique.
\end{theorem}

\qquad

\section{Proof of the main result}

In order to prove Theorem \ref{T1} we first need to prove some auxiliary
results. We start with the following lemma which is a direct adaptation of
\cite[Lemma 3.1]{julio}. This result will be used to prove Lemma \ref{L3}.

\begin{lemma}
\label{L1}Let $U\subset\mathbb{R}^{N}$ be an open $\delta$-connected set and
$u\in L^{p}\left(  U\right)  $. If%
\[
\int_{U}\int_{U}J(x,y)\left\vert u(x)-u(y)\right\vert ^{p}dydx=0\text{,}%
\]
then there exists a constant $k\in\mathbb{R}$ such that $u\left(  x\right)
=k$ $a.e.$ $x\in U.$
\end{lemma}

The next lemma will also be necessary in order to prove Lemma \ref{L3}. Lemma
\ref{L2} is crucial and presents the greatest challenge in adapting the ideas
developed in \cite{julio}.

\begin{lemma}
\label{L2}Let $1<p<\infty$ and $u_{n}:\Omega\rightarrow\mathbb{R}$ be a
sequence such that $u_{n}\rightarrow0$ strongly in $L^{p}(\Omega_{\ell})$ and
weakly in $L^{p}(\Omega_{n\ell})$. If in addition%
\begin{equation}
\lim_{n\rightarrow\infty}\int_{\Omega_{n\ell}}\int_{\Omega}J(x,y)\left\vert
u_{n}(x)-u_{n}(y)\right\vert ^{p}dydx=0, \label{eq:nolcalomega}%
\end{equation}
then%
\begin{equation}
\lim_{n\rightarrow\infty}\int_{\Omega_{n\ell}}\left\vert u_{n}(x)\right\vert
^{p}dx=0, \label{eq:solonl}%
\end{equation}
that is, $u_{n}\rightarrow0$ in $L^{p}(\Omega_{n\ell})$ and hence in
$L^{p}(\Omega)$.
\end{lemma}

\textit{Proof}. First we prove that
\[
\lim_{n\rightarrow\infty}\int_{\Omega_{n\ell}}\int_{\Omega}J(x,y)\left\vert
u_{n}(x)\right\vert ^{p}dydx=0,
\]
Let $M\in\mathbb{N}_{0}$ such that $M+1<p\leq M+2$. Then, using inequality
(III) in \cite[Page 71]{L} we get%
\[%
\begin{array}
[c]{l}%
\int_{\Omega_{n\ell}}\int_{\Omega}J(x,y)\left\vert u_{n}(x)-u_{n}%
(y)\right\vert ^{p}dydx\medskip\\
=\int_{\Omega}J(x,y)\left\vert u_{n}(x)-u_{n}(y)\right\vert ^{p-M}\left\vert
u_{n}(x)-u_{n}(y)\right\vert ^{M}dydx\medskip\\
\geq C_{p}\int_{\Omega_{n\ell}}\int_{\Omega}J(x,y)\left(  \left\vert
u_{n}(x)\right\vert ^{p-M-2}u_{n}(x)-\left\vert u_{n}(y)\right\vert
^{p-M-2}u_{n}(y)\right)  \medskip\\
\quad\times\left(  u_{n}(x)-u_{n}(y)\right)  \left\vert u_{n}(x)-u_{n}%
(y)\right\vert ^{M}dydx\medskip\\
\geq C_{p}\int_{\Omega_{n\ell}}\int_{\Omega}J(x,y)\left(  \left\vert
u_{n}(x)\right\vert ^{p-M-2}u_{n}(x)-\left\vert u_{n}(y)\right\vert
^{p-M-2}u_{n}(y)\right)  \medskip\\
\quad\times\left(  u_{n}(x)-u_{n}(y)\right)  \left(  \left\vert u_{n}%
(x)\right\vert -\left\vert u_{n}(y)\right\vert \right)  ^{M}dydx\medskip\\
=C_{p}\int_{\Omega_{n\ell}}\int_{\Omega}J(x,y)\medskip\\
\quad\times\left[  \left\vert u_{n}(x)\right\vert ^{p-M}+\left\vert
u_{n}(y)\right\vert ^{p-M}-u_{n}(x)u_{n}(y)\left(  \left\vert u_{n}%
(x)\right\vert ^{p-M-2}+\left\vert u_{n}(y)\right\vert ^{p-M-2}\right)
\right]  \medskip\\
\quad\times\sum_{j=0}^{M}\binom{M}{j}(-1)^{j}\left\vert u_{n}(y)\right\vert
^{j}\left\vert u_{n}(x)\right\vert ^{M-j}dydx\medskip\\
=C_{p}\left(  \sum_{j=0}^{M}\binom{M}{j}(-1)^{j}\int_{\Omega_{n\ell}}%
\int_{\Omega}J(x,y)\left\vert u_{n}(x)\right\vert ^{p-j}\left\vert
u_{n}(y)\right\vert ^{j}dydx\right.  \medskip\\
\quad+\sum_{j=0}^{M}\binom{M}{j}(-1)^{j}\int_{\Omega_{n\ell}}\int_{\Omega
}J(x,y)\left\vert u_{n}(y)\right\vert ^{p-M+j}\left\vert u_{n}(x)\right\vert
^{M-j}dydx\medskip\\
\quad+\sum_{j=0}^{M}\binom{M}{j}(-1)^{j+1}\int_{\Omega_{n\ell}}\int_{\Omega
}J(x,y)\left\vert u_{n}(x)\right\vert ^{p-j-2}u_{n}(x)\left\vert
u_{n}(y)\right\vert ^{j}u_{n}(y)dydx\medskip\\
\quad+\left.  \sum_{j=0}^{M}\binom{M}{j}(-1)^{j+1}\int_{\Omega_{n\ell}}%
\int_{\Omega}J(x,y)\left\vert u_{n}(y)\right\vert ^{p-M+j-2}u_{n}(y)\left\vert
u_{n}(x)\right\vert ^{M-j}u_{n}(x)dydx\right)  \medskip\\
=C_{p}\left(  \int_{\Omega_{n\ell}}\int_{\Omega}J(x,y)\left\vert
u_{n}(x)\right\vert ^{p}dydx+\sum_{j=1}^{M}\binom{M}{j}(-1)^{j}\left\langle
T_{J,p/j}(|u_{n}|^{j}),|u_{n}|^{p-j}\right\rangle \right.  \medskip\\
\quad+\int_{\Omega_{n\ell}}\int_{\Omega}J(x,y)\left\vert u_{n}(y)\right\vert
^{p}dydx+\sum_{j=0}^{M-1}\binom{M}{j}(-1)^{j}\left\langle T_{J,p/(p-M+j)}%
(|u_{n}|^{p-M+j}),|u_{n}|^{M-j}\right\rangle \medskip\\
\quad+\sum_{j=0}^{M}\binom{M}{j}(-1)^{j+1}\left\langle T_{J,p/(j+1)}%
(\left\vert u_{n}\right\vert ^{j}u_{n}),\left\vert u_{n}\right\vert
^{p-j-2}u_{n}\right\rangle \medskip\\
\quad+\left.  \sum_{j=0}^{M}\binom{M}{j}(-1)^{j+1}\left\langle
T_{J,p/(p-M+j-1)}(\left\vert u_{n}\right\vert ^{p-M+j-2}u_{n}),\left\vert
u_{n}\right\vert ^{M-j}u_{n}\right\rangle \right)  .
\end{array}
\]
Since $u_{n}\rightarrow0$ weakly in $L^{p}(\Omega_{n\ell})$, if $0<r<p$, we
get that $|u_{n}|^{r}$ is bounded in $L^{p/r}(\Omega_{n\ell})$, and then
$|u_{n}|^{r}\rightarrow0$ weakly in $L^{p/r}(\Omega_{n\ell})$. By compactness
of $T_{J,p/r}$ in $L^{p/r}(\Omega_{n\ell})$, we have that
\[
T_{J,p/r}(|u_{n}|^{r})\rightarrow0\quad\text{and}\quad T_{J,p/r}(|u_{n}%
|^{r-1}u_{n})\rightarrow0,
\]
both convergences in $L^{p/r}(\Omega_{n\ell})$. On the other hand
\[
(p/r)^{\prime}=\frac{p/r}{p/r-1}=\frac{p/r}{(p-r)/r}=p/(p-r),
\]
and since $|u_{n}|^{p-r}$ converges weakly to zero in $L^{p/(p-r)}%
(\Omega_{n\ell})$, by \cite[Proposition 3.5]{brezis}, we get
\[
\left\langle T_{J,p/r}(|u_{n}|^{r}),|u_{n}|^{p-r}\right\rangle \rightarrow
0,\qquad\left\langle T_{J,p/r}(|u_{n}|^{r-1}u_{n}),|u_{n}|^{p-r-1}%
u_{n}\right\rangle \rightarrow0.
\]
Now, we observe that $M+1<p\leq M+2$. Then: if $1\leq j\leq M$, $1<p/j$; if
$0\leq j<M$, $1<p/(p-M+j)$; if $0\leq j\leq M$, $1<p/(j+1)$ and
$1<p/(p-M+j-1)$. Therefore
\begin{align*}
\text{if }1  &  \leq j\leq M,\qquad\lim_{n\rightarrow\infty}\left\langle
T_{J,p/j}(|u_{n}|^{j}),|u_{n}|^{p-j}\right\rangle =0,\\
\text{if }0  &  \leq j<M,\qquad\lim_{n\rightarrow\infty}\left\langle
T_{J,p/(p-M+j)}(|u_{n}|^{p-M+j}),|u_{n}|^{M-j}\right\rangle =0,\\
\text{if }0  &  \leq j\leq M,\qquad\lim_{n\rightarrow\infty}\left\langle
T_{J,p/(j+1)}(\left\vert u_{n}\right\vert ^{j}u_{n}),\left\vert u_{n}%
\right\vert ^{p-j-2}u_{n}\right\rangle =0,\\
\text{if }0  &  \leq j\leq M,\qquad\lim_{n\rightarrow\infty}\left\langle
T_{J,p/(p-M+j-1)}(\left\vert u_{n}\right\vert ^{p-M+j-2}u_{n}),\left\vert
u_{n}\right\vert ^{M-j}u_{n}\right\rangle =0.
\end{align*}
Then
\[%
\begin{array}
[c]{l}%
0=\lim_{n\rightarrow\infty}\int_{\Omega_{n\ell}}\int_{\Omega}J(x,y)\left\vert
u_{n}(x)-u_{n}(y)\right\vert ^{p}dydx\medskip\\
\geq C_{p}\lim_{n\rightarrow\infty}\left(  \int_{\Omega_{n\ell}}\int_{\Omega
}J(x,y)\left\vert u_{n}(x)\right\vert ^{p}dydx+\int_{\Omega_{n\ell}}%
\int_{\Omega}J(x,y)\left\vert u_{n}(y)\right\vert ^{p}dydx\right) \\
\geq0,
\end{array}
\]
and thus
\[
\lim_{n\rightarrow\infty}\int_{\Omega_{n\ell}}\int_{\Omega}J(x,y)\left\vert
u_{n}(x)\right\vert ^{p}dydx=0.
\]
Let us next define
\[
A_{\delta}^{0}:=\left\{  x\in\Omega_{n\ell}:dist(x,\Omega_{\ell}%
)<\delta\right\}  .
\]
Notice that thanks to property (3) and to the fact that $\Omega_{n\ell}$ is
open we see that $A_{\delta}^{0}$ is open and nonempty. In particular it has
positive $N$-dimensional measure. For any $x\in\overline{A_{\delta}^{0}}$ we
consider the continuous and \textit{strictly} positive function
$g(x):=|B_{2\delta}(x)\cap\Omega_{\ell}|$. Since $\overline{A_{\delta}^{0}}$
is a compact set, there exists a constant $m>0$ such that $g(x)\geq m$ for any
$x\in\overline{A_{\delta}^{0}}$. As a consequence
\begin{align*}
\int_{\Omega_{n\ell}}\int_{\Omega_{\ell}}J(x-y)|u_{n}(x)|^{p}dydx  &  \geq
\int_{A_{\delta}^{0}}\int_{B_{2\delta}(x)\cap\Omega_{\ell}}J(x-y)|u_{n}%
(x)|^{p}dydx\\
&  \geq mC\int_{A_{\delta}^{0}}|u_{n}(x)|^{p}dx.
\end{align*}
Therefore, thanks to \eqref{eq:solonl}, $u_{n}\rightarrow0$ in $L^{p}%
(A_{\delta}^{0})$. In order to iterate this argument we notice that at this
point we know that $u_{n}\rightarrow0$ strongly in $A_{\delta}^{0}$ and weakly
in $\Omega_{n\ell}\setminus\overline{A_{\delta}^{0}}$, hence again from
\eqref{eq:nolcalomega} we get
\begin{equation}
\lim_{n\rightarrow\infty}\int_{\Omega_{n\ell}\setminus\overline{A_{\delta}%
^{0}}}\int_{A_{\delta}^{0}}J(x-y)|u_{n}(x)|^{p}dydx=0. \label{eq:soloadelta}%
\end{equation}
Since $\Omega_{n\ell}$ is $\delta$ connected, $dist(\Omega_{n\ell}%
\setminus\overline{A_{\delta}^{0}},A_{\delta}^{0})<\delta$. Considering now
\[
A_{\delta}^{1}=\{x\in\Omega_{n\ell}\setminus\overline{A_{\delta}^{0}%
}:dist(x,A_{\delta}^{0})<\delta\},
\]
and proceeding as before, we obtain, from \eqref{eq:soloadelta}, that
$u_{n}\rightarrow0$ strongly in $A_{\delta}^{1}$. This argument can be
repeated and gives strong converge in $L^{p}(A_{\delta}^{j})$ for
\[
A_{\delta}^{j}=\left\{  x\in\Omega_{n\ell}\setminus\overline{\bigcup_{0\leq
i<j}A_{\delta}^{i}}:dist\Big(x,\bigcup_{0\leq i<j}A_{\delta}^{i}%
\Big)<\delta\right\}  .
\]
Since $\Omega_{n\ell}$ is bounded, we have, for a finite number $J\in
\mathbb{N}$,
\[
\Omega_{n\ell}=\bigcup_{0\leq i<J}A_{\delta}^{i}%
\]
and therefore the proof is complete. \qed

\begin{lemma}
\label{L3}There is a constant $c>0$ such that
\[
\int_{\Omega_{\ell}}\frac{\left\vert \nabla u\right\vert ^{p}}{p}+\frac{1}%
{p}\int_{\Omega_{n\ell}}\int_{\mathbb{R}^{N}}J(x,y)\left\vert
u(x)-u(y)\right\vert ^{p}dydx\geq c\left\Vert u\right\Vert _{L^{p}\left(
\Omega\right)  }^{p}%
\]
for all $u\in W$.
\end{lemma}

\textit{Proof}. We proceed by contradiction. Assume there exists $u_{n}\in H$
such that $\left\Vert u_{n}\right\Vert _{L^{p}\left(  \Omega\right)  }=1$ and
\[
\int_{\Omega_{\ell}}\frac{\left\vert \nabla u_{n}\right\vert ^{p}}{p}+\frac
{1}{p}\int_{\Omega_{n\ell}}\int_{\mathbb{R}^{N}}J(x,y)\left\vert
u_{n}(x)-u_{n}(y)\right\vert ^{p}dydx\rightarrow0.
\]
Then, $\int_{\Omega_{\ell}}\left\vert \nabla u_{n}\right\vert ^{p}%
\rightarrow0$ and $\int_{\Omega_{n\ell}}\int_{\mathbb{R}^{N}}J(x,y)\left\vert
u_{n}(x)-u_{n}(y)\right\vert ^{p}dydx\rightarrow0$. Since $u_{n}$ is bounded
in $L^{p}\left(  \Omega_{\ell}\right)  $ and $\int_{\Omega_{\ell}}\left\vert
\nabla u_{n}\right\vert ^{p}\rightarrow0$, by the Sobolev imbedding theorem,
passing to a subsequence we get that $u_{n}\rightarrow k_{1}$ in
$W^{1,p}\left(  \Omega_{\ell}\right)  $ for some $k_{1}\in\mathbb{R}$. We
argue next in the nonlocal part $\Omega_{n\ell}$. Since $u_{n}$ is bounded in
$L^{p}\left(  \Omega_{n\ell}\right)  $, passing to another subsequence we have
that $u_{n}\rightharpoonup u$ in $L^{p}\left(  \Omega_{n\ell}\right)  $.
Furthermore, since%
\[
\int_{\Omega_{n\ell}}\int_{\mathbb{R}^{N}}J(x,y)\left\vert u_{n}%
(x)-u_{n}(y)\right\vert ^{p}dydx\rightarrow0,
\]
we get that the limit $u$ verifies that
\begin{align}
&  \int_{\Omega_{n\ell}}\int_{\Omega_{n\ell}}J(x,y)\left\vert
u(x)-u(y)\right\vert ^{p}dydx\label{1}\\
&  \leq\lim_{n\rightarrow\infty}\int_{\Omega_{n\ell}}\int_{\Omega_{n\ell}%
}J(x,y)\left\vert u_{n}(x)-u_{n}(y)\right\vert ^{p}dydx=0,\nonumber
\end{align}
and%
\begin{align}
&  \int_{\Omega_{n\ell}}\int_{\Omega_{\ell}}J(x,y)\left\vert
u(x)-u(y)\right\vert ^{p}dydx\label{2}\\
&  \leq\lim_{n\rightarrow\infty}\int_{\Omega_{n\ell}}\int_{\Omega_{\ell}%
}J(x,y)\left\vert u_{n}(x)-u_{n}(y)\right\vert ^{p}dydx=0.\nonumber
\end{align}
From (\ref{1}), using Lemma \ref{L1} and the fact that $\Omega_{n\ell}$ is an
open $\delta$-connected set, we deduce that $u=k_{2}$ in $\Omega_{n\ell}$ for
some $k_{2}\in\mathbb{R}$. On the other side, from (\ref{2}) we obtain
\[
\int_{\Omega_{n\ell}}\int_{\Omega_{\ell}}J(x,y)\left\vert k_{1}-k_{2}%
\right\vert ^{p}dydx=0
\]
and so, recalling conditions (3) and (J1) we must have $k_{1}=k_{2}$. We next
see that $k_{1}=0$. We have that $u_{n}=0$ in $\mathbb{R}^{N}\diagdown\Omega$.
If $\partial\Omega\cap\partial\Omega_{\ell}\not =\emptyset$, then
$u_{n\mid\partial\Omega\cap\partial\Omega_{\ell}}=0$; and from the convergence
$u_{n}\rightarrow k_{1}$ in $W^{1,p}\left(  \Omega_{\ell}\right)  $, we
conclude that $k_{1}=0$. If $\partial\Omega\cap\partial\Omega_{\ell}%
=\emptyset$, then in this case we have that $dist\left(  \Omega_{n\ell
},\mathbb{R}^{N}\diagdown\Omega\right)  =0$. Now, using that $u_{n}=0$ in
$\mathbb{R}^{N}\diagdown\Omega$,
\[
\int_{\Omega_{n\ell}}\int_{\mathbb{R}^{N}\diagdown\Omega}J(x,y)\left\vert
u_{n}\left(  x\right)  \right\vert ^{p}dydx\rightarrow0
\]
and $u_{n}\rightharpoonup k_{2}$ in $L^{p}\left(  \Omega_{n\ell}\right)  $, we
derive that $k_{2}=0$. Summing up, we have proved that $u_{n}\rightarrow0$ in
$W^{1,p}\left(  \Omega_{\ell}\right)  $ and $u_{n}\rightharpoonup0$ in
$L^{p}\left(  \Omega_{n\ell}\right)  $. Then, Lemma \ref{L2} says that
$u_{n}\rightarrow0$ in $L^{p}\left(  \Omega\right)  $. Since $1=\left\Vert
u_{n}\right\Vert _{L^{p}\left(  \Omega\right)  }$ for all $n$ we get a
contradiction. \qed

\qquad

We are now in position to prove the Theorem \ref{T1}

\textit{Proof of Theorem \ref{T1}}. By hypothesis $(f)$ we have
\[
\left\vert \widetilde{F}\left(  x,u\right)  \right\vert \leq\int
_{0}^{\left\vert u\right\vert }a(x)+b(x)|\xi|^{q}d\xi\leq a(x)\left\vert
u\right\vert +b(x)|u|^{q+1},
\]
where $a\in L^{s}(\Omega)$, for some $s>p^{\prime}$, $b\in L^{\infty}(\Omega)$
and $q<p-1$. 

By Lemma \ref{L3} we have that%
\begin{align*}
E\left(  u\right)   &  \geq C\left\Vert u\right\Vert _{L^{p}\left(
\Omega\right)  }^{p}-\int_{\Omega}\widetilde{F}\left(  x,u\right)  \,dx\\
&  \geq C\Vert u\Vert_{L^{p}(\Omega)}^{p}-\int_{\Omega}a(x)\left\vert
u\right\vert +b(x)|u|^{q+1}\text{ }dx\\
&  \geq C\Vert u\Vert_{L^{p}(\Omega)}^{p}-\left(  C^{\prime}\Vert
a\Vert_{L^{p^{\prime}}(\Omega)}+\Vert b\Vert_{L^{\infty}(\Omega)}\Vert
|u|^{q}\Vert_{L^{p^{\prime}}(\Omega)}\right)  \Vert u\Vert_{L^{p}(\Omega)}\\
&  \geq C\Vert u\Vert_{L^{p}(\Omega)}^{p}-C^{\prime}\Vert a\Vert
_{L^{p^{\prime}}(\Omega)}\Vert u\Vert_{L^{p}(\Omega)}-C^{\prime\prime}\Vert
u\Vert_{L^{p}(\Omega)}^{q}\Vert u\Vert_{L^{p}(\Omega)}%
\end{align*}
for some $C,C^{\prime},C^{\prime\prime}>0$, and so, since $q<p-1$, $E$ is
bounded from below and coercive. 

Let $1<r=\min\{p/s^{\prime},p/(q+1)\}$. Then $\widetilde{F}(\cdot,u)\in
L^{r}(\Omega)$ for all $u\in L^{p}(\Omega)$. Suppose now that $\{u_{n}\}$ is a
sequence that converges weakly to a function $u$ in $W$. On one hand the
functional $F$ given by
\[
F(u)=\int_{\Omega_{\ell}}\frac{\left\vert \nabla u\right\vert ^{p}}{p}%
+\frac{1}{p}\int_{\Omega_{n\ell}}\int_{\mathbb{R}^{N}}J(x,y)\left\vert
u(x)-u(y)\right\vert ^{p}dydx,
\]
is convex and therefore $F$ is weakly lower semicontinuous. On the other side,
we can take a subsequence $\{u_{n_{j}}\}$ such that
\[
-\lim_{j\rightarrow\infty}\int_{\Omega}\widetilde{F}(x,u_{n_{j}}%
)\,dx=\liminf_{n\rightarrow\infty}-\int_{\Omega}\widetilde{F}(x,u_{n})\,dx.
\]
Also, since $\{u_{n}\}$ is bounded in $L^{p}\left(  \Omega\right)  $,
$\left\{  \widetilde{F}(\cdot,u_{n})\right\}  $ is bounded in $L^{r}(\Omega)$
and so there exists some $\left\{  u_{n_{j_{k}}}\right\}  $ such that
$\left\{  \widetilde{F}(\cdot,u_{n_{j_{k}}})\right\}  $ converges weakly to
$\widetilde{F}(\cdot,u)$ in $L^{r}\left(  \Omega\right)  $. Hence,
\begin{gather*}
\liminf_{n\rightarrow\infty}E(u_{n})\geq\liminf_{n\rightarrow\infty}%
F(u_{n})+\liminf_{n\rightarrow\infty}-\int_{\Omega}\widetilde{F}%
(x,u_{n})\,dx\\
\geq\liminf_{n\rightarrow\infty}F(u_{n})-\lim_{k\rightarrow\infty}\int
_{\Omega}\widetilde{F}(x,u_{n_{j_{k}}})\,dx\geq F(u)-\int_{\Omega}%
\widetilde{F}(x,u)\,dx=E(u),
\end{gather*}
and $E$ is a weakly lower semicontinuous functional. Therefore, it is easy to
check that there exists a minimizer $u\in W$ by the direct method of the
calculus of variations. Next, we prove that $u$ is a weak solution of
(\ref{ec1}) and (\ref{ec2}). Let $\phi$ be a smooth function with $\phi=0$ in
$\mathbb{R}^{N}\diagdown\Omega$. Then, for all $t\in\mathbb{R}$ we have that
$\frac{\partial}{\partial t}E\left(  u+t\phi\right)  \left\vert _{t=0}\right.
=0$. In other words,%
\[%
\begin{array}
[c]{l}%
\int_{\Omega}f\left(  x,u\right)  \phi=\int_{\Omega_{\ell}}\left\vert \nabla
u\right\vert ^{p-2}\nabla u\nabla\phi\smallskip\\
\quad+\int_{\Omega_{n\ell}}\int_{\mathbb{R}^{N}}J(x,y)\left\vert
u(x)-u(y)\right\vert ^{p-2}\left(  u(x)-u(y)\right)  \left(  \phi
(x)-\phi(y)\right)  dydx.
\end{array}
\]
Now, we observe that%
\[%
\begin{array}
[c]{l}%
\int_{\Omega_{n\ell}}\int_{\mathbb{R}^{N}}J(x,y)\left\vert
u(x)-u(y)\right\vert ^{p-2}\left(  u(x)-u(y)\right)  \left(  \phi
(x)-\phi(y)\right)  dydx\smallskip\\
=\int_{\Omega_{n\ell}}\int_{\Omega_{n\ell}}J(x,y)\left\vert
u(x)-u(y)\right\vert ^{p-2}\left(  u(x)-u(y)\right)  \left(  \phi
(x)-\phi(y)\right)  dydx\smallskip\\
\quad+\int_{\Omega_{n\ell}}\int_{\mathbb{R}^{N}\diagdown\Omega_{n\ell}%
}J(x,y)\left\vert u(x)-u(y)\right\vert ^{p-2}\left(  u(x)-u(y)\right)  \left(
\phi(x)-\phi(y)\right)  dydx,
\end{array}
\]
and so, using that $J$ is symmetric and Fubini's theorem we get%
\[%
\begin{array}
[c]{l}%
\int_{\Omega_{n\ell}}\int_{\Omega_{n\ell}}J(x,y)\left\vert
u(x)-u(y)\right\vert ^{p-2}\left(  u(x)-u(y)\right)  \left(  \phi
(x)-\phi(y)\right)  dydx\smallskip\\
=-2\int_{\Omega_{n\ell}}\int_{\Omega_{n\ell}}J(x,y)\left\vert
u(x)-u(y)\right\vert ^{p-2}\left(  u(y)-u(x)\right)  dy\text{ }\phi\left(
x\right)  dx.
\end{array}
\]
On the other side,
\[%
\begin{array}
[c]{l}%
\int_{\Omega_{n\ell}}\int_{\mathbb{R}^{N}\diagdown\Omega_{n\ell}%
}J(x,y)\left\vert u(x)-u(y)\right\vert ^{p-2}\left(  u(x)-u(y)\right)  \left(
\phi(x)-\phi(y)\right)  dydx\smallskip\\
=-\int_{\Omega_{n\ell}}\int_{\mathbb{R}^{N}\diagdown\Omega_{n\ell}%
}J(x,y)\left\vert u(x)-u(y)\right\vert ^{p-2}\left(  u(y)-u(x)\right)
dy\text{ }\phi\left(  x\right)  dx\smallskip\\
\quad-\int_{\mathbb{R}^{N}\diagdown\Omega_{n\ell}}\int_{\Omega_{n\ell}%
}J(x,y)\left\vert u(x)-u(y)\right\vert ^{p-2}\left(  u(y)-u(x)\right)
dy\text{ }\phi\left(  x\right)  dx.
\end{array}
\]
Therefore, recalling that $u=\phi=0$ in $\mathbb{R}^{N}\diagdown\Omega$ we
have that%
\[%
\begin{array}
[c]{l}%
\int_{\Omega}f\left(  x,u\right)  \phi=\int_{\Omega_{\ell}}\left\vert \nabla
u\right\vert ^{p-2}\nabla u\nabla\phi\smallskip\\
\quad-2\int_{\Omega_{n\ell}}\int_{\Omega_{n\ell}}J(x,y)\left\vert
u(x)-u(y)\right\vert ^{p-2}\left(  u(y)-u(x)\right)  dy\text{ }\phi\left(
x\right)  dx\smallskip\\
\quad-\int_{\Omega_{n\ell}}\int_{\mathbb{R}^{N}\diagdown\Omega_{n\ell}%
}J(x,y)\left\vert u(x)-u(y)\right\vert ^{p-2}\left(  u(y)-u(x)\right)
dy\text{ }\phi\left(  x\right)  dx\smallskip\\
\quad-\int_{\Omega_{\ell}}\int_{\Omega_{n\ell}}J(x,y)\left\vert
u(x)-u(y)\right\vert ^{p-2}\left(  u(y)-u(x)\right)  dy\text{ }\phi\left(
x\right)  dx,
\end{array}
\]
and then $u$ is a weak solution of (\ref{ec1}) and (\ref{ec2}). Finally if
$\xi\rightarrow f\left(  x,\xi\right)  $ is strictly concave for $a.e.$
$x\in\Omega$, then $E$ is a strictly convex functional in $W$ and the
minimizer $u$ is unique. \qed
\smallskip

\textit{Acknowledgement}. We would like to thank to Julio Rossi for suggesting
us this problem.


\begin{thebibliography}{99}                                                                                               %


\bibitem {julio}G. Acosta, F. Bersetche, J. Rossi, \textit{Local and nonlocal
energy-based coupling models}, SIAM J. Math. Anal. \textbf{54} (2022), 6288--6322.

\bibitem {ABR21}G. Acosta, F. Bersetche, J. Rossi, \textit{Coupling local and
nonlocal equations with Neumann boundary conditions}, Rev. Uni\'{o}n Mat.
Argent. \textbf{65} (2023), 533--565.

\bibitem {ABR23}G. Acosta, F. Bersetche, J. Rossi, \textit{A domain
decomposition scheme for couplings between local and nonlocal equations},
Comput. Methods Appl. Math. \textbf{23} (2023), 817--830.

\bibitem {BCRR15}H. Berestycki, A. Coulon, J. Roquejoffre, L. Rossi,
\textit{The effect of a line with nonlocal diffusion on Fisher-KPP
propagation}, Math. Models Methods Appl. Sci. \textbf{25} (2015), 2519--2562.

\bibitem {brezis}H. Brezis, \textit{Functional analysis, Sobolev spaces and
partial differential equations}. Universitext. Springer, New York (2011).

\bibitem {CR23}M. Capanna, J. Rossi, \textit{Mixing local and nonlocal
evolution equations.} Mediterr. J. Math. \textbf{20}, 59 (2023).

\bibitem {CCR06}E. Chasseigne, M. Chaves, J. Rossi, \textit{Asymptotic
behavior for nonlocal diffusion equations}, J. Math. Pures Appl. \textbf{86}
(2006), 271--291.

\bibitem {CY}J. Chen, Y. Tang, \textit{Homogenization of non-local nonlinear
p-Laplacian equation with variable index and periodic structure}, J. Math.
Phys. \textbf{64}., 061502 (2023).

\bibitem {CERW07}C. Cortazar, M. Elgueta, J. Rossi, N. Wolanski,
\textit{Boundary fluxes for nonlocal diffusion}, J Differ. Equations
\textbf{234} (2007), 360--390.

\bibitem {CERW08}C. Cortazar, M. Elgueta, J. Rossi, N. Wolanski, \textit{How
to approximate the heat equation with Neumann boundary conditions by nonlocal
diffusion problems}, Arch. Ration. Mech. Anal. \textbf{187} (2008), 137--156.

\bibitem {CD07}J. Coville, L. Dupaigne, \textit{On a non-local equation
arising in population dynamics}, Proc. R. Soc. Edinb., Sect. A, Math.
\textbf{137} (2007), 727--755.

\bibitem {DB}M. D'Elia, P. Bochev, \textit{Formulation, analysis and
computation of an optimization-based local-to-nonlocal coupling method},
Results Appl. Math. \textbf{9} (2021), 100--129.

\bibitem {DLLT}Q. Du, X. Li, J. Lu, X. Tian, \textit{A quasi-nonlocal coupling
method for nonlocal and local diffusion models.} SIAM J. Numer. Anal. 56
(2018), no. 3, 1386-1404.

\bibitem {ERS}C. Esteve, J. Rossi, A. San Antolin, \textit{Upper bounds for
the decay rate in a nonlocal p-Laplacian evolution problem}, Boundary Value
Probl. \textbf{2014 }(2014), 1--10.

\bibitem {GW}C. Gal, M. Warma, \textit{Nonlocal transmission problems with
fractional diffusion and boundary conditions on non-smooth interfaces.}
Commun. Partial Differ Equations \textbf{42} (2017), 579--625.

\bibitem {GQR}A. G\'{a}rriz, F. Quir\'{o}s, J. Rossi, \textit{Coupling local
and nonlocal evolution equations}, Calc. Var. Partial Differ. Equ. \textbf{59}
(2020), 1--24.

\bibitem {GT}B. Gess, J. T\"{o}lle, \textit{Ergodicity and local limits for
stochastic local and nonlocal p-Laplace equations}, SIAM J. Math. Anal.
\textbf{48} (2016), 4094--4125.

\bibitem {GS}G. Graham, I. Sloan, \textit{On the Compactness of Certain
Integral Operators}, J. Math Anal. Appl. \textbf{68 }(1979), 580--594.

\bibitem {IR}L. Ignat, J. Rossi, \textit{Decay estimates for nonlocal problems
via energy methods}, J. Math. Pures Appl. \textbf{92} (2009), 163--187.

\bibitem {IRA}L. Ignat, J. Rossi, A. Antolin, \textit{Decay estimates for
nonlinear nonlocal diffusion problems in the whole space}, J. Anal. Math.
\textbf{122} (2014), 375--401.

\bibitem {K15}D. Kriventsov, \textit{Regularity for a Local-Nonlocal
Transmission Problem}, Arch. Ration. Mech. Anal. \textbf{217 }(2015), 1103--1195.

\bibitem {L}P. Lindqvist, \textit{Notes on the p-Laplace equation}. Report.
University of Jyv\"{a}skyl\"{a}. Department of Mathematics and Statistics 102.
(ISBN 951-39-2586-2/pbk). 80 p. (2006).

\bibitem {MRT}J. Maz\'{o}n, J. Rossi, J. Toledo. \textit{Fractional
p-Laplacian evolution equations, }J. Math. Pures Appl. \textbf{105} (2016), 810--844.

\bibitem {MD14}T. Mengesha, Q. Du, \textit{The bond-based peridynamic system
with Dirichlet-type volume constraint}, Proc. R. Soc. Edinb., Sect. A, Math.
\textbf{144} (2014), 161--186.

\bibitem {SOR}B. dos Santos, S. Oliva, J. Rossi, \textit{A local/nonlocal
diffusion model}, Appl. Anal. \textbf{101} (2022), 5213--5246.

\bibitem {SS}M. Schuster, V. Schulz, \textit{Interface Identification
constrained by Local-to-Nonlocal Coupling.} arXiv preprint arXiv:2402.12871 (2024).

\bibitem {SL10}S. Silling, R. Lehoucq, \textit{Peridynamic theory of solid
mechanics.} Adv. Appl. Mech. \textbf{44} (2010), 73--168.

\bibitem {SLL}J. Sun, J. Li, Q. Liu, \textit{Cauchy problem of a nonlocal
p-Laplacian evolution equation with nonlocal convection}, Nonlinear Anal.,
Theory Methods Appl., Ser. A, Theory Methods \textbf{95} (2014), 691--702.

\bibitem {W02}X. Wang, \textit{Metastability and stability of patterns in a
convolution model for phase transitions.} J. Differ. Equations \textbf{183}
(2002), 434--461.

\bibitem {Z04}L. Zhang, \textit{Existence, uniqueness and exponential
stability of traveling wave solutions of some integral differential equations
arising from neuronal networks}, J. Differ. Equations \textbf{197} (2004), 162--196.
\end{thebibliography}
\end{document}